\documentclass[smallcondensed]{svjour3}                     
\smartqed  
\usepackage{graphicx}
\usepackage{booktabs}
\usepackage{amssymb} 
\usepackage{amsfonts}
\usepackage{url}
\usepackage[leqno]{amsmath}
\usepackage{epsfig}
\usepackage[numbers,sort]{natbib}
\usepackage{tabularx}
\usepackage{comment}
\usepackage{mathtools}
\usepackage{multirow}
\usepackage{breakcites}
\usepackage{geometry}
\usepackage{lscape}
\newcommand{\be}{\begin{equation}}
\newcommand{\ee}{\end{equation}}

\newcommand{\alg}[1]{\texttt{#1}}
\newcommand{\vx}{\mathbf{x}}
\newcommand{\myvec}[1]{\mathbf{#1}}
\newcommand{\vy}{\mathbf{y}}
\newcommand{\eg}{e.g.}
\newcommand{\ie}{i.e.}

\newcommand{\ouralg}{\alg{\textsc{WOO}}}

\newcommand\myeq{\mathrel{\stackrel{\makebox[0pt]{\mbox{\normalfont\tiny def}}}{=}}}

\newcommand{\IUnaryAddInd}{$I^{1}_{\epsilon+}(\mathcal{Y}^t_*)$}

\newcommand{\fct}[1]{\mathbf{f}(#1)}
\newcommand{\norm}[1]{||#1||}
\newcommand{\cheb}{Tchebycheff~}
\newcommand{\gfct}{g_{\mathbf{w}}}
\spnewtheorem{Assumption}{Assumption}{\bf}{\it}

\usepackage{algorithm}
\usepackage{algpseudocode}
\begin{document}

\title{Revisiting Norm Optimization for Multi-Objective Black-Box Problems: A Finite-Time Analysis}

\titlerunning{Revisiting Norm Optimization for Multi-Objective Problems}        

\author{Abdullah Al-Dujaili\and S. Suresh}


\institute{Abdullah Al-Dujaili \at
              CSAIL, MIT, Cambridge, MA 02139, USA\\
              \email{aldujail@mit.edu}           
           \and
          S. Suresh \at
                            School of Computer Engineering, Nanyang Technological University, Singapore 639798\\ 
                            \email{ssundaram@ntu.edu.sg}
}

\date{Received: date / Accepted: date}

\maketitle

\begin{abstract}
The complexity of Pareto fronts  imposes a great challenge on the convergence analysis of multi-objective optimization methods. While most theoretical convergence studies have addressed finite-set and/or discrete problems, others have provided probabilistic guarantees, assumed a total order on the solutions, or studied their asymptotic behaviour. In this paper, we revisit the Tchebycheff weighted method in a hierarchical bandits setting and provide  a finite-time bound on the Pareto-compliant additive $\epsilon$-indicator. To the best of our knowledge, this paper is one of few that establish a link between weighted sum methods and quality indicators in finite time.
\end{abstract}

\noindent\textbf{Problem.} This paper is concerned with the {Multi-Objective Black-Box} Optimization (\textsc{MOBBO}) problem given a {finite number of function evaluations}. With $n$ decision variables and $m$ objectives, the problem  has the mathematical form:

\begin{equation}
\begin{aligned}
& {\text{minimize}}
& & \myvec{y} = \myvec{f}(\myvec{x}) = (f_1(\myvec{x}), \ldots, f_m(\myvec{x}))\\
& \text{where}
& & \myvec{x} =  (x_1, \ldots, x_n)\in \mathcal{X} \subset \mathbb{R}^n\\
& & & \myvec{y} =  (y_1, \ldots, y_m)\in \mathcal{Y} \subset \mathbb{R}^m\;,
\end{aligned}
\label{eq:problem_def}
\end{equation}
where $\myvec{x}$ is called the \emph{decision vector (solution)}, $\myvec{y}$ is called the \emph{objective vector}, $\mathcal{X}$ 
is the \emph{feasible decision space}
, and $\mathcal{Y}={\text{ $\bigtimes$}}_{1 \leq j\leq m} \mathcal{Y}_j$ is the reachable \emph{objective space}, where $\mathcal{Y}_j$ is the $j$th-objective space. By black-box, we mean to say that there is no closed-form expression of $\mathbf{f}$ and that its derivatives are neither symbolically nor numerically available. However, $\myvec{f}$ can be evaluated point-wise, but each evaluation is typically expensive in terms of computational resources (\eg, time, power, money).

In practice, the objective functions $\{f_j\}_{1\leq j \leq m}$ are conflictual. Therefore, the problem may have  a set of incomparable optimal solutions: each is inferior to the other in some
objectives and superior in other objectives, inducing a partial order on the set of
feasible solutions.


\vspace{0.5em}
\noindent\textbf{Related Work.} The notion of Pareto-optimality was first introduced in the field of engineering in 1970~\cite{stadler1979survey}. The \textsc{MOBBO} solvers can be broadly classified into generative or preference methods based on the decision maker's role. The former does not require any inputs from the decision maker in solving the problem (only selects solution at the end), whereas preference methods require input from decision maker at the beginning.  The knowledge of decision maker may affect the solution quality in the preference method with regard to some objectives, but it reduces the complexity in \textsc{MOBBO} solver. Conventionally, most commonly used \textsc{MOBBO} solvers convert multiple objectives into a single (or a series of) objective optimization problem \cite{zadeh1963optimality,messac2004normal}.  The adaptive scalarization~\cite{kim2005adaptive} overcomes the problem of non-uniform distribution of optimal solution and the non-convex Pareto fronts in scalarizing approaches~\cite{das1997closer}. 

With regard to black-box optimization, optimistic algorithms---whose foundations come from the multi-armed bandit theory~\cite{Munos2011,al2016naive,Al-Dujaili2016,Al-Dujaili2016a}---consider partitions of the search space at multiple scales in search of optimal solutions. These methods enjoy provable finite-time performance and asymptotic convergence. On the other hand, research works with the \cheb metric have lacked theoretical guarantees~\cite{van2014multi}. To this end, this paper aims to bridge the gap in understanding the theoretical underpinnings of the weighted \cheb method and link the convergence of scalarization methods to Pareto-compliant quality indicators in a multi-arm bandits setting.

\vspace{0.5em}
\noindent\textbf{Our Contributions.}  This paper addresses a class of weighted sum methods for MOBBO problems.
While most of the literature work has established the asymptotic optimality of weighted sum methods to a single Pareto-optimal solution under certain conditions, our theoretical contributions here are of two-fold. First, we show that the weighted \cheb problem for Lipschitz MOBBO is as well Lipschitz. Second, we present a finite-time upper bound on the Pareto-compliant quality indicator of the approximation set obtained from solving the weighted \cheb problem capturing its convergence to the whole Pareto front. All of this is motivated by the success of the optimism in the face of uncertainty principle that helps us employ an optimistic method, which we refer to as the Weighted Optimistic Optimization (\ouralg) algorithm.  The sequential decision-making approach in \ouralg~formulates the weighted decision space $\mathcal{X}$ as a hierarchy of simple bandit problems over subspaces of  $\mathcal{X}$ and looks for the optimal solution through $\mathcal{X}$-partitioning search trees. At each step, \ouralg~expands the subspace which may contain the optimum. Based on smoothness assumption, the convergence analysis is presented in finite time and validated empirically using a set of synthetic problems.

\vspace{0.5em}
\noindent\textbf{Paper Organization.} The rest of the paper is organized as follows. First, a brief introduction to basic concepts in multi-objective optimization is provided, along with the notations and terminology used through out the paper. Then, we motivate treating the weighted \cheb problem in a multi-arm bandits setting by proving its smoothness. This is followed by introducing the weighted optimistic optimization algorithm that exploits the smoothness without the need for its knowledge. Furthermore, theoretical and empirical analysis of the proposed algorithm is presented. Towards the end, we conclude with a discussion on potential future research investigations.

\section{Formal Background}
\label{sec:formal_bg}

This section presents basic concepts in multi-objective optimizations. First, the notion of Pareto dominance is described. Second, approaches to assess the performance of multi-objective solvers are discussed.  Third, we formally define the weighted \cheb problem.

\subsection{Pareto Dominance}
\label{sec:pareto-dom}

An objective vector $\myvec{y}^1\in \mathcal{Y}$ is more preferable than another vector $\myvec{y}^2\in \mathcal{Y}$, if $\myvec{y}^1$ is at least as good as $\myvec{y}^2$ in all objectives \emph{and} better with respect to at least one objective. $\myvec{y}^1$ is then said to be \emph{dominating} $\myvec{y}^2$. This notion of dominance is commonly known as \textit{Pareto dominance}~\cite{pareto-book}, which leads to a \textit{partial order} on the objective space, where we can define a Pareto optimal vector to be one that is non-dominated by any other vector in $\mathcal{Y}$.  Nevertheless,  $\myvec{y}^1$ and $\myvec{y}^2$ may be incomparable to each other, because each is inferior to the other in some objectives and superior in other objectives. Hence, there can be several Pareto optimal vectors. This concept is presented in the following definitions~\cite{Zitzler2003,Loshchilov2013}.

\begin{definition}[Pareto dominance] The vector $\myvec{y}^1$ dominates the vector $\myvec{y}^2$, that is to say, $\myvec{y}^1 \prec\myvec{y}^2$$\iff$ $y^1_j\leq y^2_j$ for all $j\in \{1,\ldots,m\}$ and $y^1_k < y^2_k$ for at least one $k\in\{1,\ldots,m\}$.
  \label{def:parteo_dominance}
\end{definition}

\begin{definition}[Strict Pareto dominance] The vector $\myvec{y}^1$ strictly dominates the vector~$\myvec{y}^2$ if $\myvec{y}^1$ is better than $\myvec{y}^2$ in all the objectives, that is to say, $\myvec{y}^1 \prec\prec\myvec{y}^2\iff$ $y^1_j< y^2_j$ for all $j\in \{1,\ldots,m\}$.
  \label{def:strict_parteo_dominance}
\end{definition}

\begin{definition}[Weak Pareto dominance] The vector $\myvec{y}^1$ weakly dominates the vector~$\myvec{y}^2$ if $\myvec{y}^1$ is not worse than $\myvec{y}^2$ in all the objectives, that is to say, $\myvec{y}^1 \preceq\myvec{y}^2\iff$ $y^1_j\leq y^2_j$ for all $j\in \{1,\ldots,m\}$.
  \label{def:weak_parteo_dominance}
\end{definition}

\begin{definition}[Pareto optimality of vectors]
  \label{def:paretoptimal}
  Let $\hat{\myvec{y}}\in \mathcal{Y}$ be a vector. $\hat{\myvec{y}}$ is Pareto optimal $\iff$ $\nexists \myvec{y} \in \mathcal{Y}$ such that  $\myvec{y}\prec \hat{\myvec{y}}$. The set
  of all Pareto optimal vectors is referred to as the Pareto front and denoted as $\mathcal{Y}^*$. The corresponding decision vectors (solutions) are referred to as the Pareto optimal solutions or the Pareto set and denoted by $\mathcal{X}^*$.
\end{definition}
Thus, the solution to the \textsc{MOBBO} problem \eqref{eq:problem_def} is its Pareto optimal solutions (Pareto front in the objective space). Practically, \textsc{\textsc{MOBBO}} solvers aim to identify a set of objective vectors that represent the Pareto front (or a good approximation of it). We refer to this set as the \emph{approximation set}. 
\begin{definition}[Approximation set]
  Let $A \subseteq \mathcal{Y}$ be
  a set of objective vectors. $A$ is called an approximation
  set if any element of $A$ does not dominate or is
  not equal to any other objective vector in $A$. The set
  of all approximation sets is denoted as $\Omega$. Note that $\mathcal{Y}^*\in \Omega$.
  \label{def:approx_set}
\end{definition}
Furthermore, denote the \emph{ideal point (utopian vector)} by $\myvec{y}^*\myeq$ $(\min_{\myvec{y}\in\mathcal{Y}^*} y_1,$
$\ldots,\min_{\myvec{y}\in\mathcal{Y}^*} y_m)$. Likewise, let us denote the (or one of the) global optimizer(s) of the $j$th objective function by $\myvec{x}^*_j$, \ie, $y^*_j=f_j(\myvec{x}^*_j)$. Note that $\myvec{x}^*_j\in \mathcal{X}^*$. 
Without loss of generality, we assume that $\myvec{y}^*$ is the zero vector.

\subsection{Performance Assessment for Multi-Objective Optimization Methods}
\label{sec:bg:bbmo}

Given two approximation sets $A,B\in \Omega$, it is not that easy to tell which set is better, particularly if their elements are incomparable~\cite{Zitzler2003}. In general, two aspects are considered in an approximation set: i). its distance (the closer the better) to
the optimal Pareto front and ii). its diversity (the higher the better) within the
optimal Pareto front. To this end, several \emph{quality indicators} have been proposed~\cite{Knowles2006}.
The quality of an approximation set is measured by a so-called (unary) quality indicator $I:\Omega\to \mathbb{R}$, assessing a specific property of the approximation set. Likewise, an $l$-ary quality indicator $I:\Omega^l\to \mathbb{R}$ quantifies quality differences between $l$ approximation sets~\cite{Zitzler2003,Custodio2011}. A quality indicator is not Pareto-compliant if it contradicts the order induced by the Pareto-dominance relations described in Section~\ref{sec:formal_bg}. One  commonly-used quality indicators is the Pareto-compliant additive $\epsilon$-indicator, which is defined formally next. 

\begin{definition}
  (Additive $\epsilon$-indicator~\cite{Zitzler2003}) For any two approximation sets $A, B \in \Omega$, the additive $\epsilon$-indicator $I_{\epsilon+}$ is defined as:
  \begin{equation}
  I_{\epsilon+}(A,B) = \inf_{\epsilon\in \mathbb{R}}\{\forall \myvec{y}^2 \in B,\; \exists \myvec{y}^1\in A : \myvec{y}^1 \preceq_{\epsilon+} \myvec{y}^2\}\,\label{eq:def_binary_epsilon}
  \end{equation}
  where $\myvec{y}^1 \preceq_{\epsilon+} \myvec{y}^2\iff y^1_j \leq \epsilon + y^2_j$ for all $j\in\{1,\ldots, m\}$. If $B$ is the Pareto front $\mathcal{Y}^*$ (or a good---in terms of diversity and closeness to the Pareto front---approximation reference set $R\in\Omega$ if $\mathcal{Y}^*$ is unknown) then $I_{\epsilon+}(A,B)$ is referred to as the unary additive  epsilon indicator and is denoted by $I^{1}_{\epsilon+}(A)$, \ie, $I^{1}_{\epsilon+}(A)\myeq I_{\epsilon+}(A,\mathcal{Y}^*)$.
  \label{def:epsilon_indicator}
\end{definition}

In essence, $I^{1}_{\epsilon+}(A)$~measures the smallest amount $\epsilon$ needed to translate each element in the Pareto front $\mathcal{Y}^*$ such that it is \emph{weakly} dominated by at least one element in the approximation set $A$. 
Note that $I^{1}_{\epsilon+}(A)\myeq I_{\epsilon+}(A,\mathcal{Y}^*)\geq0$ as no element in $A$ strictly dominates any element in $\mathcal{Y}^*$. Thus, the closer $I^{1}_{\epsilon+}(A)$ to $0$, the better the quality of $A$.

\subsection{Weighted Sum Methods}
\label{sec:weighted-methods}
In weighted sum methods, the idea is to assign a non-negative weight value $w$ for each objective and  minimize the weighted sum of the objectives. Denote the \emph{element-wise product} or \emph{Hadamard product} of two vectors $\mathbf{a}$ and $\mathbf{b}$ of the same dimensionality by $\mathbf{a} \odot \mathbf{b}$, and the $l_p$-norm of a vector $\mathbf{a}$ by $\norm{\mathbf{a}}_p$, then 
weighted sum methods have  the mathematical form $\norm{ \mathbf{w} \odot |\fct{\vx} - \mathbf{z}^*|}_p$. As the focus of this paper is the weighted \cheb problem, we define it formally next.

\begin{definition}(Weighted \cheb
  problem) Let $\mathbf{w} \in \mathbb{R}^m$ be a non-negative vector, and $\mathbf{z}^* \in \mathbb{R}^m$ be a reference point. Then, the weighted \cheb formulation of problem~\eqref{eq:problem_def} is defined as:
  \begin{equation}
  \begin{aligned}
  & {\text{minimize}}
  & &  \gfct(\vx) = \max_{1\leq j \leq m} w_j|f_j(\myvec{x}) - z^*_j| \\ 
  &\; &\; &\;\;\hspace{2em} = \norm{ \mathbf{w} \odot |\fct{\vx} - \mathbf{z}^*|}_\infty\\
  & \text{where}
  & & \myvec{x} =  (x_1, \ldots, x_n)\in \mathcal{X} \subset \mathbb{R}^n\;.
  \end{aligned}
  \label{eq:cheb_def}
  \end{equation}
  \label{def:cheb_def}
\end{definition}

Under certain conditions, the solution to problem (3) corresponds to a Pareto-optimal solution of problem (1). This is stated in the following theorem.
\begin{theorem}\cite[Theorem 3.4.5]{miettinen1999nonlinear}.
  Let $\vx$ be a Pareto-optimal solution of problem~\eqref{eq:problem_def}, then there exists a positive
  weighting vector $\mathbf{w}$ such that $\vx$ is a solution of the weighted
  \cheb problem~\eqref{eq:cheb_def} where the reference point is the
  utopian objective vector $\vy^*$.
  \label{thm:pareto-cheb}
\end{theorem}

\noindent\textbf{Proof.} See \cite[Page 98]{miettinen1999nonlinear}. $\hfill\blacksquare$
\vspace{1em}

Similar to problem~\eqref{eq:problem_def}, the weighted \cheb problem~\eqref{eq:cheb_def} is black-box but with a single objective. Therefore, using a computational budget of $v(t)$ function evaluations, we would like to devise an algorithm that searches the decision space $\mathcal{X}$ over $t$ iterations in a \emph{sequential decision making framework}, where each sample $\vx^k$ may depend on the previous sampled points and their corresponding function values $\{(\vx^i, \gfct(\vx^i))\}_{1\leq i \leq k-1}$. After the final iteration $t$ of the algorithm, we can obtain the sampled point with the best possible $\gfct$ value:
\begin{equation}
\vx(t) \in \arg\min_{1\leq i \leq v(t)} \gfct(\vx^i)\;.
\label{eq:xt}
\end{equation}
Also, with regard to the corresponding multi-objective problem at hand, we can obtain an approximation set to the Pareto front as the set of the non-dominated sampled points:
\begin{equation}
\begin{aligned}
 \mathcal{Y}^t_*=\{\mathbf{f}(\mathbf{x}^i )\; |\;   i = 1,\ldots, v(t)\;,\;\mathbf{f}(\mathbf{x}^k ) \nprec \mathbf{f}(\mathbf{x}^i ), \\ \forall  k = 1,\ldots, v(t), k\neq i \} \in \Omega\;.
 \label{eq:approx-set}
 \end{aligned}
\end{equation}

\section{The Weighted \cheb Problem: A Multi-Armed Bandit View}

In the previous section, we formulated the weighted \cheb problem as searching the decision space in a sequential decision making framework, where the present sample may depend on the previous observed samples. One popular principle of such framework is the \emph{optimism in the face of uncertainty} principle, which suggests following the optimal strategy
with respect to the most favorable scenario among all possible scenarios that are compatible
with the obtained observations about the problem at hand~\cite{Munos2014}. It has been the main principle in multi-armed bandits settings and was later extended to many (including infinite) arms under probabilistic or structural (smoothness) assumptions about the arm rewards.
In this section, we approach problem~\eqref{eq:cheb_def} in a multi-armed bandits setting assuming a structural smoothness. Our smoothness assumption is motivated by the next theorem, which shows that the \cheb function is Lipschitz-continuous if the corresponding objectives are.


\begin{theorem}
  The \cheb function for an optimization problem with $m$ Lipschitz-continuous objectives whose constants are $\{L_j\}_{1\leq j \leq m}$, respectively,  is also Lipschitz-continuous with the Lipschitz constant being less than or equal to $\sqrt{m}\cdot \max_{j}w_j L_j$.
  \label{thm:lipschitz-cheb}
\end{theorem}

\noindent\textbf{Proof}. As the objectives of the problem  are Lipschitz continuous, for any $\vx^1$ and  $\vx^2 \in \mathcal{X}$, we have 
$$|f_j(\vx^1)-f_j(\vx^2)| \leq L_j \norm{\vx^1 - \vx^2}_2\;,\; j=1,\ldots,m\;. $$
Squaring both sides of the above inequality can be expressed in vector notation as follows.
\begin{eqnarray}
\norm{\fct{\vx^1}-\fct{\vx^2}}^2_2& =& \sum_{j=1}^{m} {L^2_j} \norm{\vx^1-\vx^2}^2_2 \nonumber \\
&\leq &m\cdot \max_j L^2_{j} \norm{\vx^1-\vx^2}^2_2\;. \nonumber
\end{eqnarray}
Since the $l_\infty$-norm of any vector is bounded by its $l_2$-norm, it follows that 
\begin{eqnarray}
\norm{\fct{\vx^1}-\fct{\vx^2}}_\infty& \leq & \sqrt{m}\cdot \max_j L_j \norm{\vx^1-\vx^2}_2\;. \label{eq:linf-bound}
\end{eqnarray}
Thus, the absolute difference in the weighted \cheb function $\gfct(\vx)$ values at any two vectors $\vx^1$ and $\vx^2 \in \mathcal{X}$ for problem~\eqref{eq:problem_def} can be bounded as below.
\begin{eqnarray}
\begin{aligned}
|\gfct(\vx^1)-\gfct(\vx^2)| = & {\;}| \norm{\mathbf{w} \odot |\fct{\vx^1} - \mathbf{z}^*|}_\infty \\
& { } - \norm{\mathbf{w} \odot |\fct{\vx^2} - \mathbf{z}^*|}_\infty
| \nonumber\\
\leq & {\;}\norm{\mathbf{w} \odot |\fct{\vx^1} - \mathbf{z}^*| \\
& {} - \mathbf{w} \odot |\fct{\vx^2} - \mathbf{z}^*|}_\infty\;.
\nonumber 
\end{aligned}
\end{eqnarray}
To guarantee Pareto-optimal solutions to problem~\eqref{eq:problem_def}, Theorem~\ref{thm:pareto-cheb} suggests that $\mathbf{z}^*$ should be the utopian vector $\vy^*$, which is---without loss of generality---the zero vector based on Section~\ref{sec:pareto-dom}. Therefore, we have
\begin{eqnarray}
|\gfct(\vx^1)-\gfct(\vx^2)| 
& \leq & \norm{\mathbf{w} \odot (\fct{\vx^1}-\fct{\vx^2})}_\infty \nonumber\\
& \leq &\norm{\mathbf{w}}_\infty \norm{\fct{\vx^1}-\fct{\vx^2}}_\infty 
\nonumber \\
& \leq & \sqrt{m}\cdot \max_j w_j L_j \norm{\vx^1-\vx^2}_2\;. \nonumber
\end{eqnarray}
Therefore, $\gfct$'s Lipschitz constant is less than or equal $\sqrt{m}\cdot \max_j w_jL_j$. $\hfill\blacksquare$
\vspace{1em}

Theorem~\ref{thm:lipschitz-cheb} established the Lipschitz continuity of the weighted \cheb problem provided that all objectives are Lipschitz-continuous. This nicely supports our smoothness assumption and justifies following the optimism in the face of uncertainty principle.

Nevertheless, the aforementioned theorem does not tell us about the quality of problem~\eqref{def:cheb_def}'s solution~$\vx(t)$ and its relation to problem~\eqref{eq:problem_def}'s approximation set~$\mathcal{Y}^t_*$. Given $v(t)$ function evaluations, our next theorem upper bounds the quality of the approximation set $\mathcal{Y}^t_*$ obtained from all points sampled in $t$ iterations to solve problem~\eqref{eq:cheb_def} by the best value $\gfct$ found of these samples. 

\begin{figure}[h!]
  \begin{center}
    \renewcommand{\arraystretch}{1.2}
    \resizebox{0.9\textwidth}{!}{
    \includegraphics[scale=1]{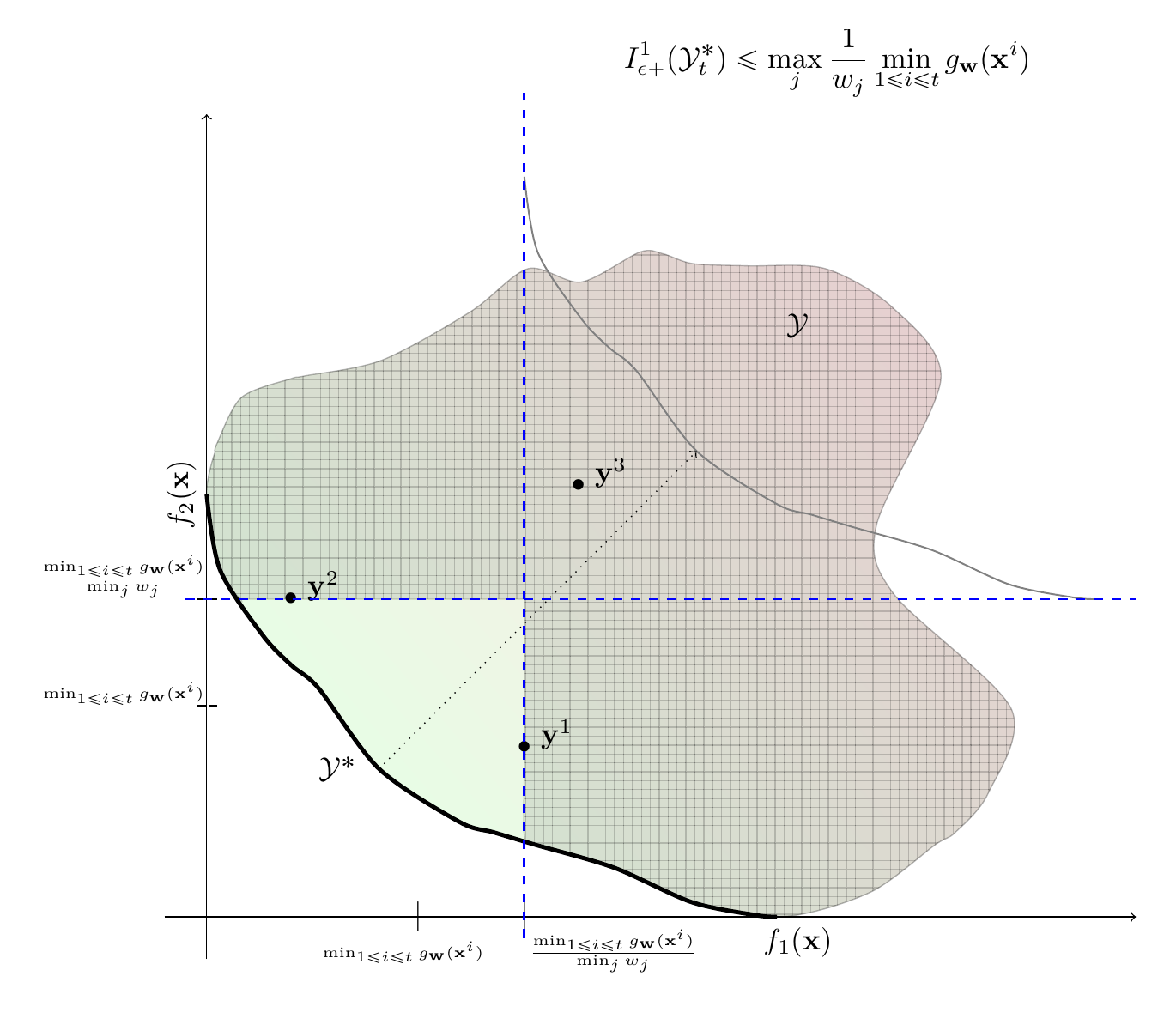}}
  \end{center}
  \caption{A visual illustration of Theorem~\ref{thm:ind-bound} in the particular case of $\min_j w_j < 1$ and $v(t)=3$ function evaluations. There exists at least one objective vector $\in \mathcal{Y}^t_*$ that lies along the top or right side of $\mathcal{Y}$'s unshaded area (e.g., $\vy^1$ and $\vy^2$).}
  \label{fig:ind-cheb-bound}    
\end{figure}

\begin{theorem}
  \IUnaryAddInd$~\leq~ \max_j \frac{1}{w_j} \gfct({\vx(t)})$.
  \label{thm:ind-bound}
\end{theorem}
\noindent\textbf{Proof}. From Definition~\ref{def:cheb_def} and~\eqref{eq:xt}, the $l_\infty$-norms of the weighted objective vectors $\{||\mathbf{w} \odot \vy^i||_\infty\}_{1\leq i \leq v(t)}$ of the sampled points $\{\vx^i\}_{1\leq i \leq v(t)}$ are greater than or equal $\gfct(\vx(t))$ (note that the ideal point $\vy^*$ is the zero vector). This suggest that there exists one vector $\vy^1 \in \mathcal{Y}^t_*$ such that for all vectors $\vy^2 \in \mathcal{Y}^*$, $$\vy^1 \preceq_{\gfct(\vx(t))/\min_j w_j} \vy^2\;.$$ Then, it follows from Definition~\ref{def:epsilon_indicator} that \IUnaryAddInd$~\leq~ \max_j \frac{1}{w_j} \gfct({\vx(t)})$. See Figure~\ref{fig:ind-cheb-bound} for pictorial proof of the above with a bi-objective problem.$\hfill\blacksquare$
\vspace{1em}

While the former theorem justifies the smoothness assumption, the latter provides a link to analyze~$\mathcal{Y}^t_*$ with respect to $\vx(t)$
The next section---motivated by the two theoretical insights at hand---presents an optimistic algorithm

\section{The Weighted Optimistic Optimization Algorithm}

\begin{figure}[ht]
  \begin{center}
    \renewcommand{\arraystretch}{1.2}
    \includegraphics[scale=1]{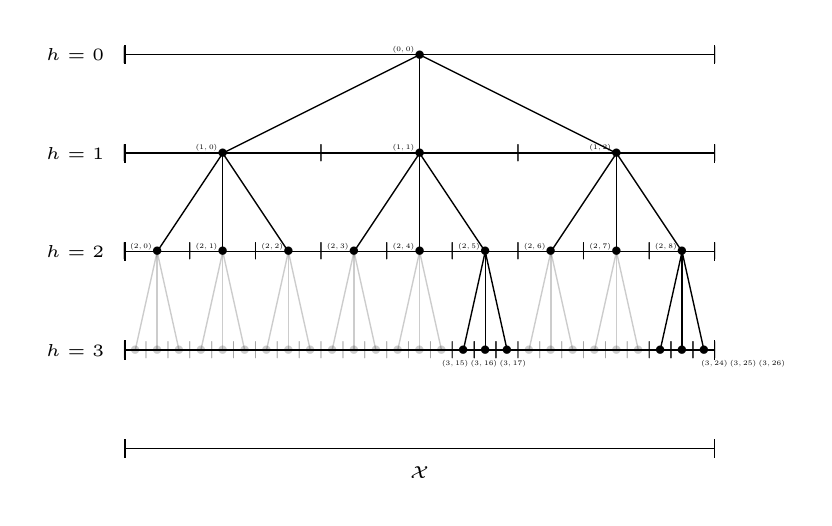}
  \end{center}
  \caption{Hierarchical bandits partitioning the decision space $\mathcal{X}$ with a partition factor of $P=3$ in a form of a $P$-ary tree, where each node is indexed by its depth $h$ and its order $i$ along the same: ($h$,$i$).} 
  \label{fig:hierarchical_partition}    
\end{figure}

Single-objective continuous optimization problems (such as~\eqref{eq:cheb_def}) can be represented as a structured bandit problem where the objective value is a function of some arm parameters. To cope with the infinitely many arms (points in $\mathcal{X}$), arms can be generated in a hierarchical fashion transforming the problem from a many-arm bandit to a hierarchy of multi-armed bandits. As shown in Figure~\ref{fig:hierarchical_partition}, one can use a space-partitioning procedure to iteratively construct finer and finer partitions of the search space $\mathcal{X}$ at multiple depths (scales)~$h\in \mathbb{N}_0$. Formally, at depth $h\geq0$, $\mathcal{X}$ can be partitioned into a set of $P^h$ cells/subspaces $\mathcal{X}_{h,i}$ where $0\leq i \leq P^{h}-1$ such that $\bigcup_{0 \leq i \leq P^h-1}\mathcal{X}_{h,i} = \mathcal{X}$. These cells are represented by nodes of a $P$-ary tree $\mathcal{T}$, where a node ($h,i$) represents the cell $\mathcal{X}_{h,i}$---the root node ($0,0$) represents the entire search space $\mathcal{X}_{0,0}=\mathcal{X}$. A parent node possesses $P$ child nodes $\{(h + 1, i_p)\}_{1\leq p \leq P}$, whose cells form a partition of the parent's cell $\mathcal{X}_{h,i}$. The set of leaves in $\mathcal{T}$ is denoted as $\mathcal{L}$. Likewise, the set of leaves at depth $h$ are denoted by $\mathcal{L}_{h}$. Each cell is has a representative point (state) $\vx_{h,i} \in \mathcal{X}_{h,i}$ at which the function is evaluated (observed). The sampled values $\gfct(\vx_{h,i})$ are employed optimistically to guide the tree partitioning.

We have proved that $\gfct$ is Lipschitz-continuous for Lipschitz multi-objective problems. Note that this smoothness is unknown in general for black-box problem.  For this, optimistic methods come to rescue. Inspired from \cite{Munos2011}, the pseudo-code of the proposed Weighted Optimistic Optimization (\ouralg)~is shown in Algorithm \ref{alg:eh} \ouralg~grows a tree $\mathcal{T}$ over $\mathcal{X}$ by expanding at most one leaf node per depth in an iterative sweep across $\mathcal{T}$'s depths/levels.  At depth $h \geq 0$, a leaf node ($h,i$) with lowest function value is expanded by partitioning its subspace along one dimension of $\mathcal{X}$. The algorithm inputs/parameters are listed as follows: i) weighted \cheb function; ii) evaluation budget $v$; iii) partition factor $P$; and iv) partitioning tree depth $h_{max}(t)$. The last two parameters contribute to exploration-vs.-exploitation trade-off.

\begin{figure*}[h!]
	\centering
  \resizebox{0.85\textwidth}{!}{
    \begin{tabular}{|cc|}
      \toprule
      \multicolumn{2}{|c|}{\textbf{After 3 function evaluations}}  \\
      \textbf{Decision Space} & \textbf{Objective Space} \\

      \includegraphics[scale=0.75,trim={2cm 0 2cm 1cm},clip]{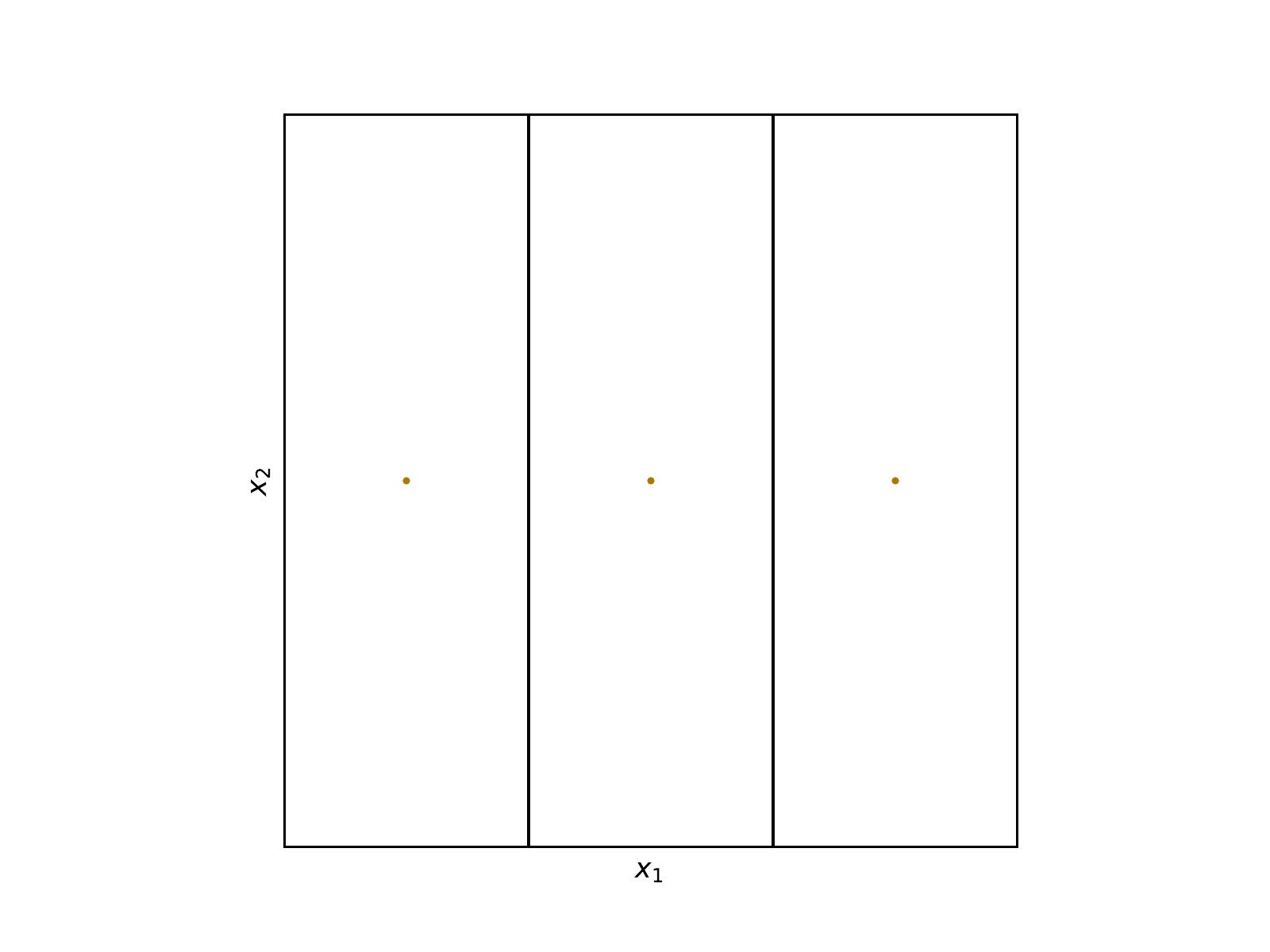} & \includegraphics[scale=0.75,trim={0 0 1cm 1cm},clip]{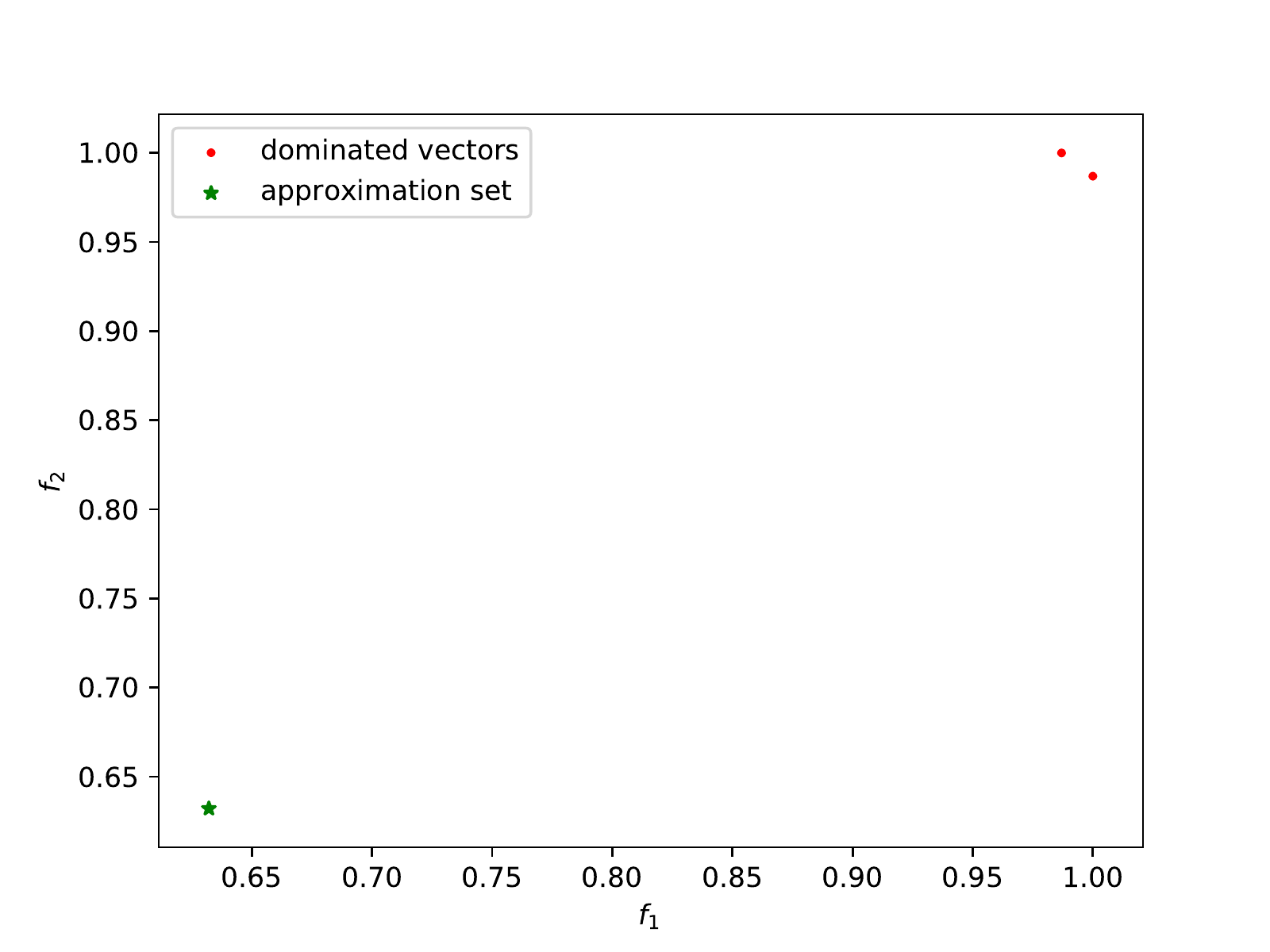} \\ \midrule
      \multicolumn{2}{|c|}{\textbf{After 10 function evaluations}} \\
      \textbf{Decision Space} & \textbf{Objective Space}\\
      \includegraphics[scale=0.75,trim={2cm 0 2cm 1cm},clip]{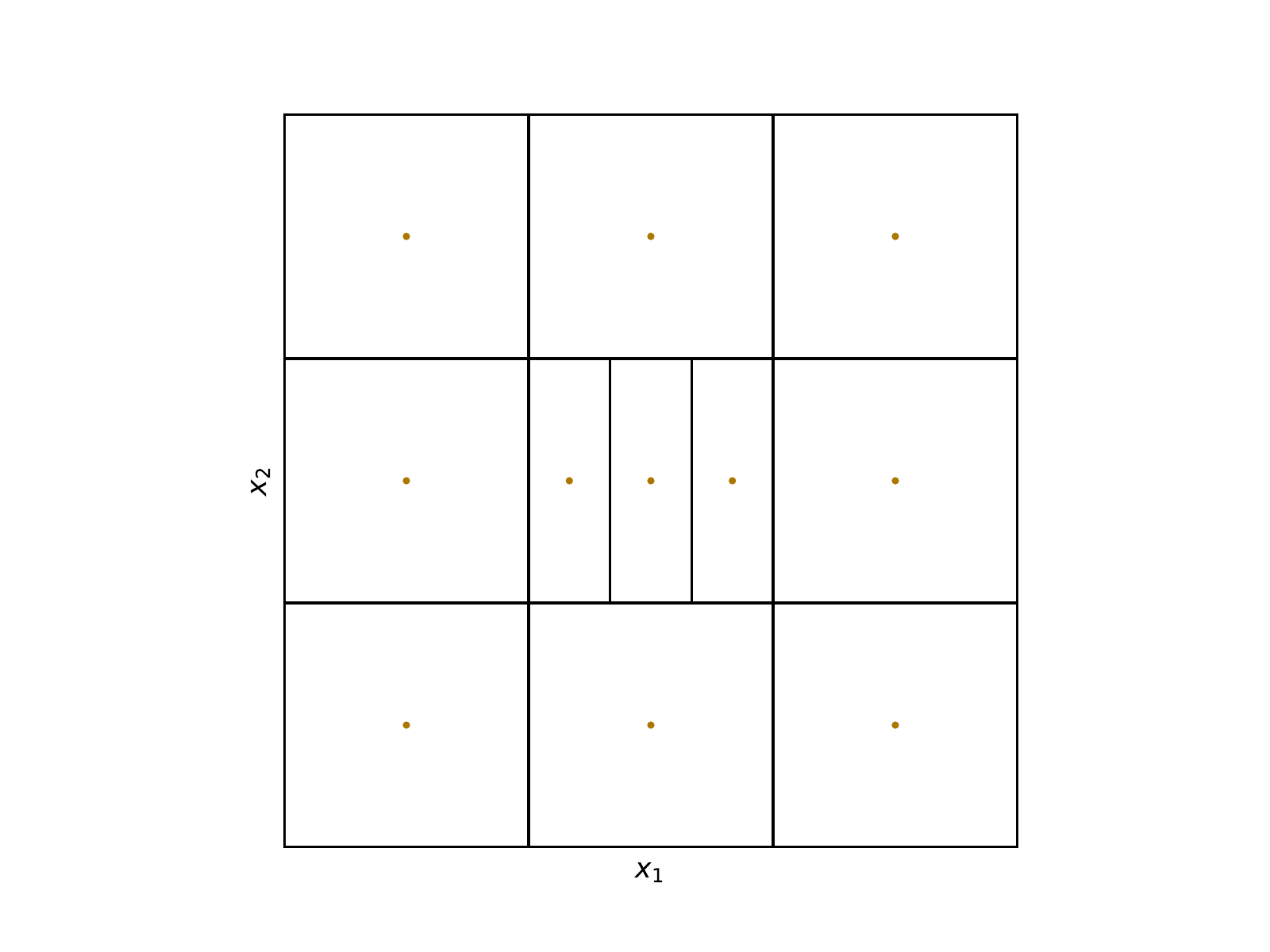} & \includegraphics[scale=0.75,trim={0 0 1cm 1cm},clip]{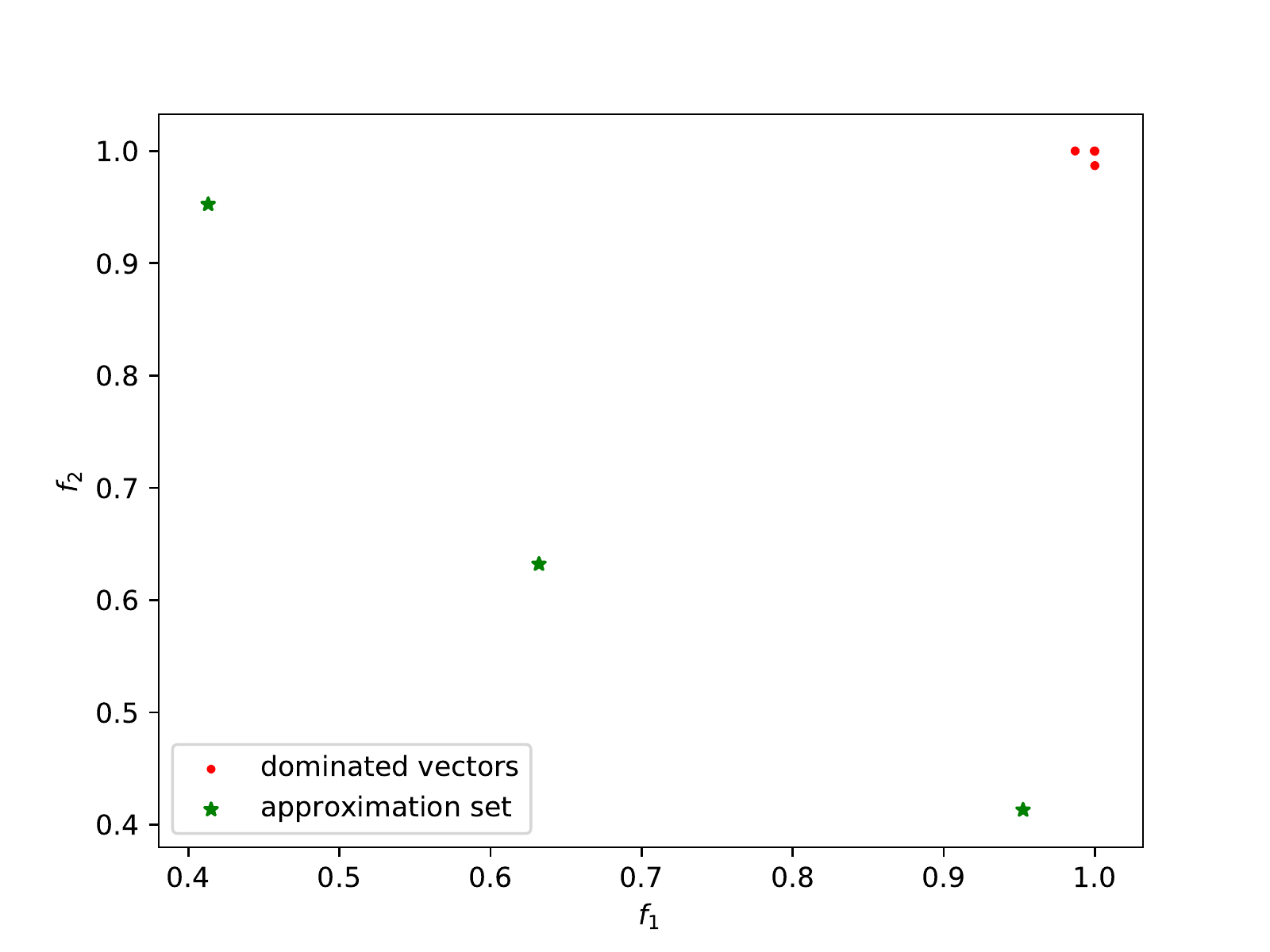}\\
      \midrule
      \multicolumn{2}{|c|}{\textbf{After 20 function evaluations}}  \\
      \textbf{Decision Space} & \textbf{Objective Space} \\
      \includegraphics[scale=0.75,trim={2cm 0 2cm 1cm},clip]{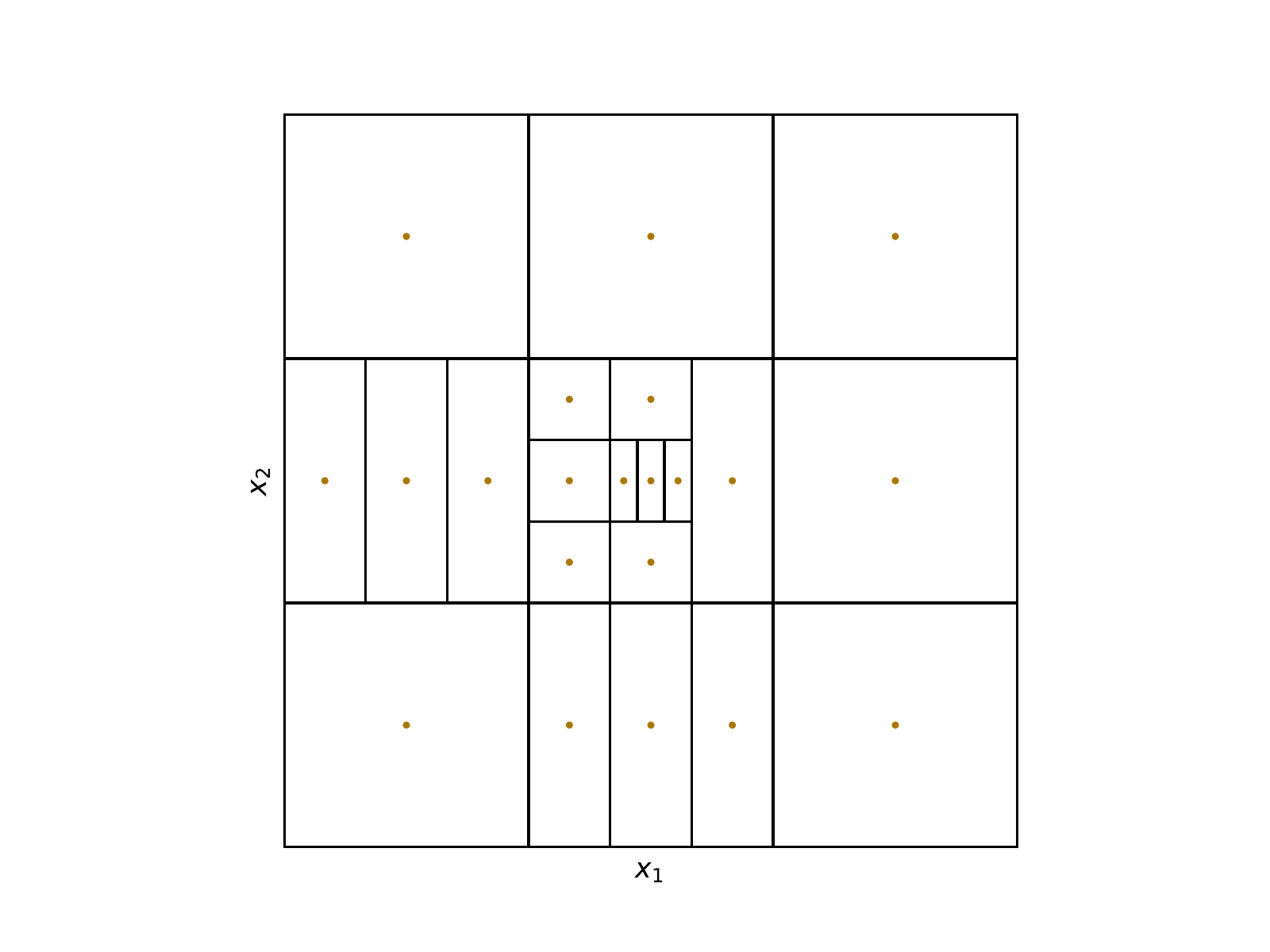} & \includegraphics[scale=0.75,trim={0 0 1cm 1cm},clip]{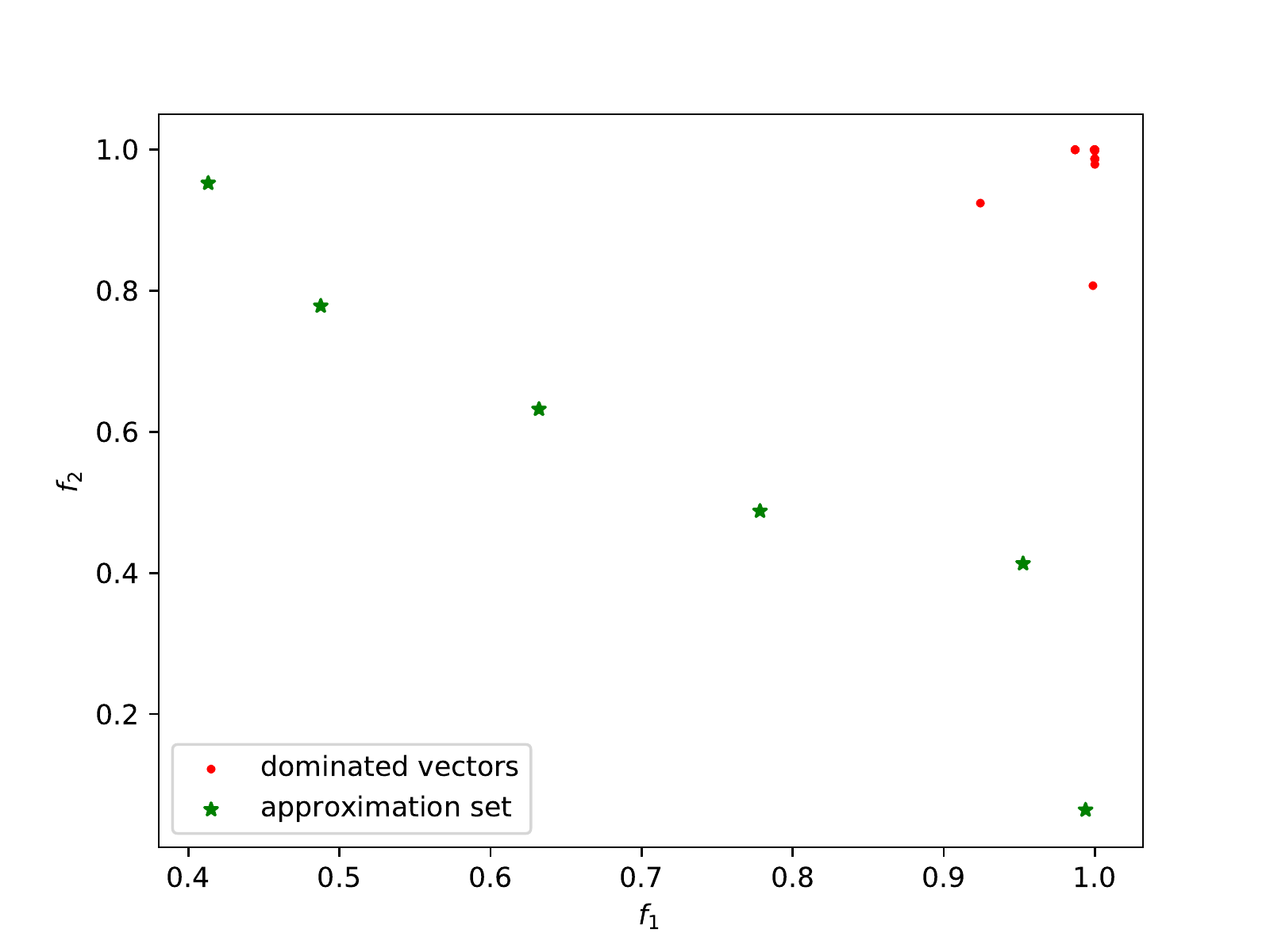} \\ \midrule
                  \multicolumn{2}{|c|}{\textbf{After 200 function evaluations}} \\
      \textbf{Decision Space} & \textbf{Objective Space}\\
      \includegraphics[scale=0.75,trim={2cm 0 2cm 1cm},clip]{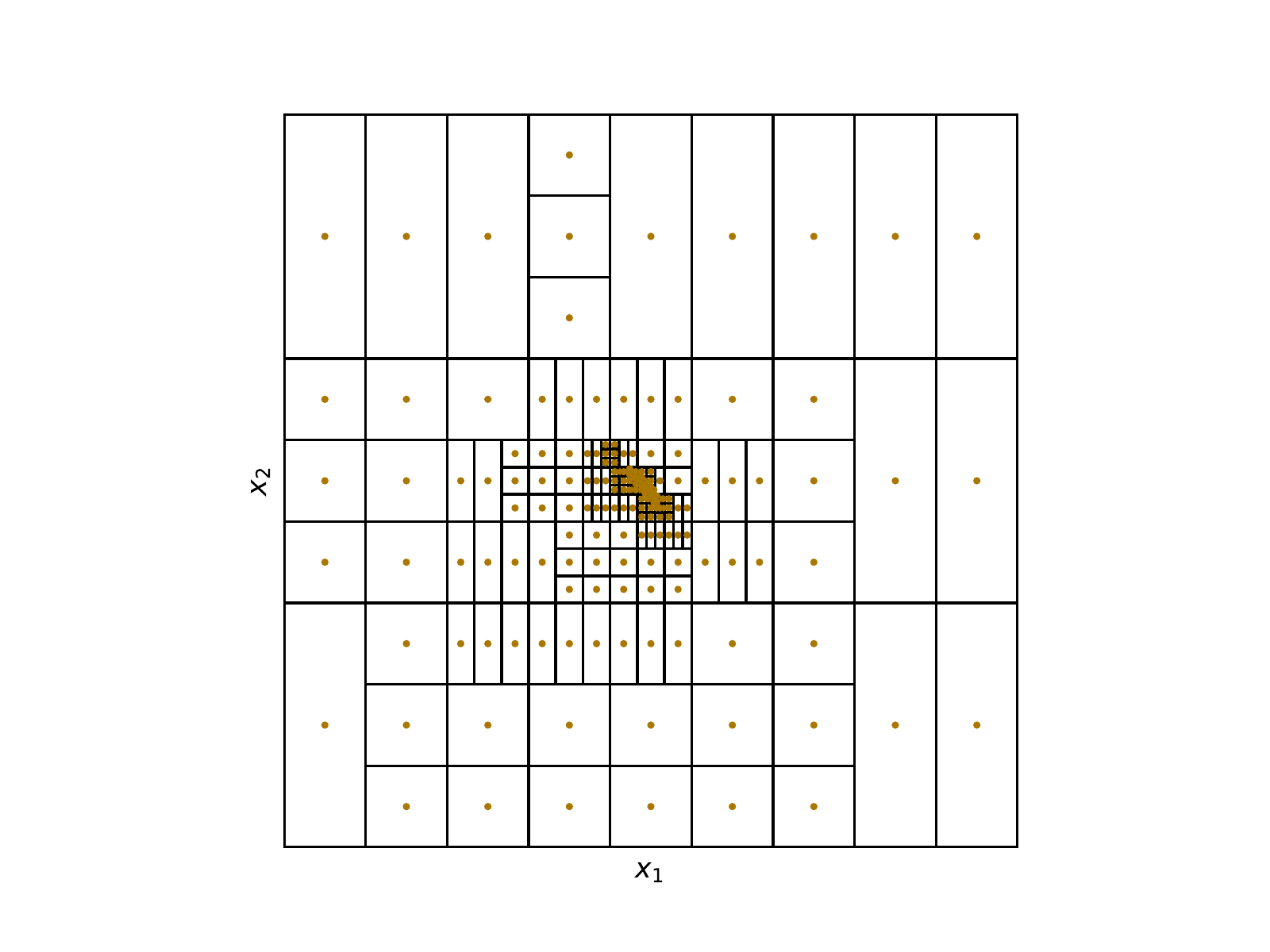} & \includegraphics[scale=0.75,trim={0 0 1cm 1cm},clip]{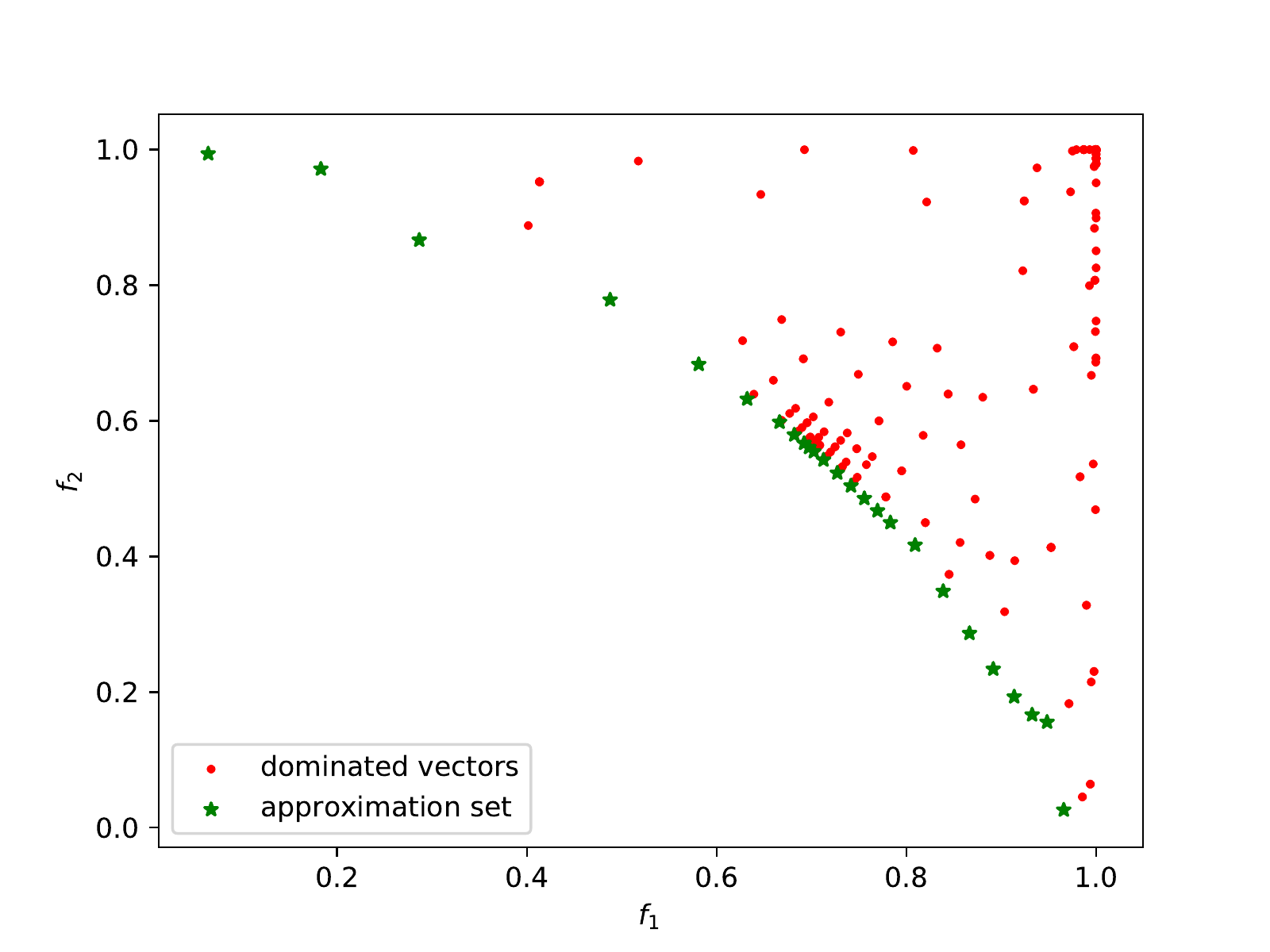}\\
      \bottomrule
    \end{tabular}
  }
  \caption{A visual illustration of \ouralg's sampling and partitioning of the decision space   on the \emph{Fonseca and Fleming} benchmark problem~\cite{fonseca1995overview}, where $\mathcal{X}=[-4,4]^2$, $f_1(\vx)=1-e^{-\sum_{i=1}^{2}(x_i - \frac{1}{\sqrt{2}})^2}$, and $f_2(\vx)=1-e^{-\sum_{i=1}^{2}(x_i + \frac{1}{\sqrt{2}})^2}$. The hierarchical partitioning is performed coordinate-wise in a round-robin fashion across iterations.}
  \label{fig:demo}
\end{figure*}

For better understanding of \ouralg, a bi-objective optimization problem is selected and Figure~\ref{fig:demo} shows working of \ouralg~at different stages and approximation of Pareto-front.

\begin{algorithm}[h!]
  \newcommand{\Statey}{\State}
  \caption{The Weighted Optimistic Optimization (\ouralg) Algorithm}  \label{alg:eh}
  \begin{algorithmic}[1]
    \Statex \textbf{Input}: 
    \Statex \hspace{1em} weighted \cheb function $\gfct$, 
    \Statex \hspace{1em} search space $\mathcal{X}$, 
    \Statex \hspace{1em} partition factor $P$,
    \Statex \hspace{1em} evaluation budget~$v$,
    \Statex \hspace{1em} partitioning tree depth $h_{max}(t)$.
    \Statex \textbf{Initialization}: 
    \Statex \hspace{1em} $t\gets 1$,
    $\mathcal{T}_1=\{(0,0)\}$,
    Evaluate $\gfct(\vx_{0, 0})$.
    \While{evaluation budget is not exhausted}
    \Statey  $\nu_{\min} \gets \infty$
    \For{$l=0$ \textbf{~to~}$\min\{\mbox{depth}{(\mathcal{T}_t}), h_{\max}(t)\}$}
    \Statey Select $(l, o)=\arg\min_{(h,i) \in \mathcal{L}_{t,l}}g^*_{h,i}$ \label{line:soo_qt}
    \If {$g^*_{l,o} < \nu_{\min}$}
    \Statey $\nu_{\min} \gets g^*_{l,o}$
    \Statey Expand $(l,o)$ into its $P$ child nodes
    \Statey Evaluate $(l,o)$'s $P$ child nodes by $\gfct$\label{ln:eval}
    \Statey Add $(l,o)$'s child nodes to $\mathcal{T}_t$
    \EndIf
    \Statey $\mathcal{T}_{t+1}\gets \mathcal{T}_t$
    \Statey $t\gets t + 1$ \label{ln:end_iteration}
    \EndFor
    \EndWhile
    \Statey \textbf{return} $\mathcal{Y}^t_*=\{\mathbf{f}(\mathbf{x}^i )\; |\;   i = 1,\ldots, v(t)\;,\;\mathbf{f}(\mathbf{x}^k ) \nprec \mathbf{f}(\mathbf{x}^i )\;, \forall  k = 1,\ldots, v, k\neq i \}$
  \end{algorithmic}
\end{algorithm}

\section{Convergence Analysis}
\label{sec:analysis}
In general, the performance of multi-armed bandit algorithms (strategies) is assessed through the notion of \emph{regret/ loss}: the difference between the present strategy's outcome and that of the optimal strategy. As our problem of interest is multi-objective, we define a Pareto-compliant regret in terms of the additive $\epsilon$-quality indicator in terms of the obtained approximation set~\eqref{eq:approx-set} as a function of the number of iterations as follows.
\begin{equation}
r(t) = I^1_{\epsilon+}(\mathcal{Y}^t_*) - I^1_{\epsilon+}(\{\mathbf{f}(\arg\min_{\vx\in \mathcal{X}} \gfct(\vx))   \}) 
\label{eq:indicator-loss}
\end{equation}

In the light of \eqref{eq:indicator-loss}, the performance of \ouralg~is analyzed two folds. First, a theoretical finite-time upper bound on the loss~\eqref{eq:indicator-loss} is proved. Second, numerical experiments are setup to validate the proven performance on a set of synthetic bi-objective problems.

\subsection{Bounding the Indicator Loss}

Building on Theorem~\ref{thm:lipschitz-cheb}, we employ the smoothness assumptions used in~\cite{Munos2011,Al-Dujaili2016}. Subsequently, \ouralg's exploration of the decision space can be quantified in terms of the near-optimality dimension $d$~\cite[Definition 1]{Munos2011}. From that, the following theorem is deduced.

\begin{theorem}
  Given $t$ iterations, let us write $h(t)$ the smallest integer $h$ such that
  \begin{equation}
  Ch_{max}(t)
  \sum_{l=0}^{h}\delta(l)^{-d} \geq t\;,
  \label{eq:loss-condition}
  \end{equation}
  where $d$ is the near-optimality dimension and $\delta(l)$ is a decreasing sequence in $l>0$ that captures $\gfct$ smoothness over the hierarchical partitioning $\{\mathcal{X}_{h,i}\}_{h\geq 0, 0\leq i < P^h}$
  Then the loss is bounded as
  \begin{equation} 
  r(t)\leq \max_j \frac{1}{w_j}\delta(
  \min(h(t),h_{max}(t)+1))\;.
  \label{eq:loss-bound}
  \end{equation}
  \label{thm:finite-time-loss}
\end{theorem}
\noindent\textbf{Proof.} Based on \cite[Theorem 2]{Munos2011}, the loss in the quality of $\gfct(\vx(t))$ is bounded by $\delta(h(t), h_{max}(t) + 1)$, \ie,
$$\gfct(\vx(t)) - \min_{\mathcal{X}} \gfct(\vx) \leq \delta(h(t), h_{max}(t) + 1)\;.$$
Multiplying both sides by $\max_j \frac{1}{w_j}$ yields
\begin{equation}
\begin{aligned}
\underbrace{\max_j \frac{1}{w_j}\gfct(\vx(t))}_{\geq I^1_{\epsilon+}(\mathcal{Y}^t_*)} - \underbrace{\max_j \frac{1}{w_j}\min_{\mathcal{X}} \gfct(\vx)}_{\geq I^1_{\epsilon+}(\{\mathbf{f}(\arg\min_{\vx\in \mathcal{X}} \gfct(\vx))  }
\\  \leq \max_j \frac{1}{w_j}\delta(h(t), h_{max}(t) + 1)\;,\nonumber
\end{aligned}
\end{equation}
Also, based on Theorem~\ref{thm:ind-bound}, the two terms on the left side are lower bounded by terms under the braces, respectively. Therefore, the loss $r(t)$ of Eq.~\eqref{eq:indicator-loss} can be bounded as follows.
\begin{equation}
\begin{aligned}
r(t) = I^1_{\epsilon+}(\mathcal{Y}^t_*) - I^1_{\epsilon+}(\{\mathbf{f}(\arg\min_{\vx\in \mathcal{X}} \gfct(\vx))   \}) 
\\ \leq \max_j \frac{1}{w_j}\delta(h(t), h_{max}(t) + 1)\;.\nonumber
\end{aligned}
\end{equation}
~$\hfill\blacksquare$

\subsection{Numerical Validation}

In this section, the presented finite-time loss bound in Theorem~\ref{thm:finite-time-loss} is validated numerically. To this end, we used the setup described in~\cite{al2016multi}: a set of synthetic bi-objective ($m=2$) problems of the form
\begin{equation}
\begin{aligned}
f_1(\myvec{x}) & = & ||\myvec{x}-\vx^*_1||^{\alpha_1}_\infty\;, \\
f_2(\myvec{x})&=&||\myvec{x}-\vx^*_2||^{\alpha_2}_\infty\;, 
\end{aligned}
\label{eq:synthetic-problem}
\end{equation}
where $\mathcal{X}=[-1,1]^n$, $n\in \{1,2\}$, and $(\vx^*_1, \vx^*_2) \in \{
(\mathbf{0},\mathbf{1}), (\mathbf{0.21}, \mathbf{0.81}), (\mathbf{0.47}, \mathbf{0.61}), (\mathbf{0.57}, \mathbf{0.57})
\}$ to capture a range of conflicting objectives. The theoretical bounds~\eqref{eq:loss-condition} and \eqref{eq:loss-bound} for the problem instances at hand were coded in Python using the \texttt{SymPy} package~\cite{meurer2017sympy} and computed as described in~\cite{al2016multi}. 

Moreover, the numerical indicator values $I^1_{\epsilon+}(\mathcal{Y}^t_*)$ (first term of ~\eqref{eq:indicator-loss}) at each iteration $t$  are computed using a Python implementation of the \ouralg~algorithm using a budget of $v=10^3$ function evaluations. On the other hand, the second term of~\eqref{eq:indicator-loss}, $I^1_{\epsilon+}(\{\mathbf{f}(\arg\min_{\vx\in \mathcal{X}} \gfct(\vx))   \})$, is computed for each of the eight problem instances using a budget of $v=5\times 10^5$ function evaluations.

\begin{landscape}
	\begin{figure}[htpb]
		\resizebox{\textwidth}{!}{
			\includegraphics{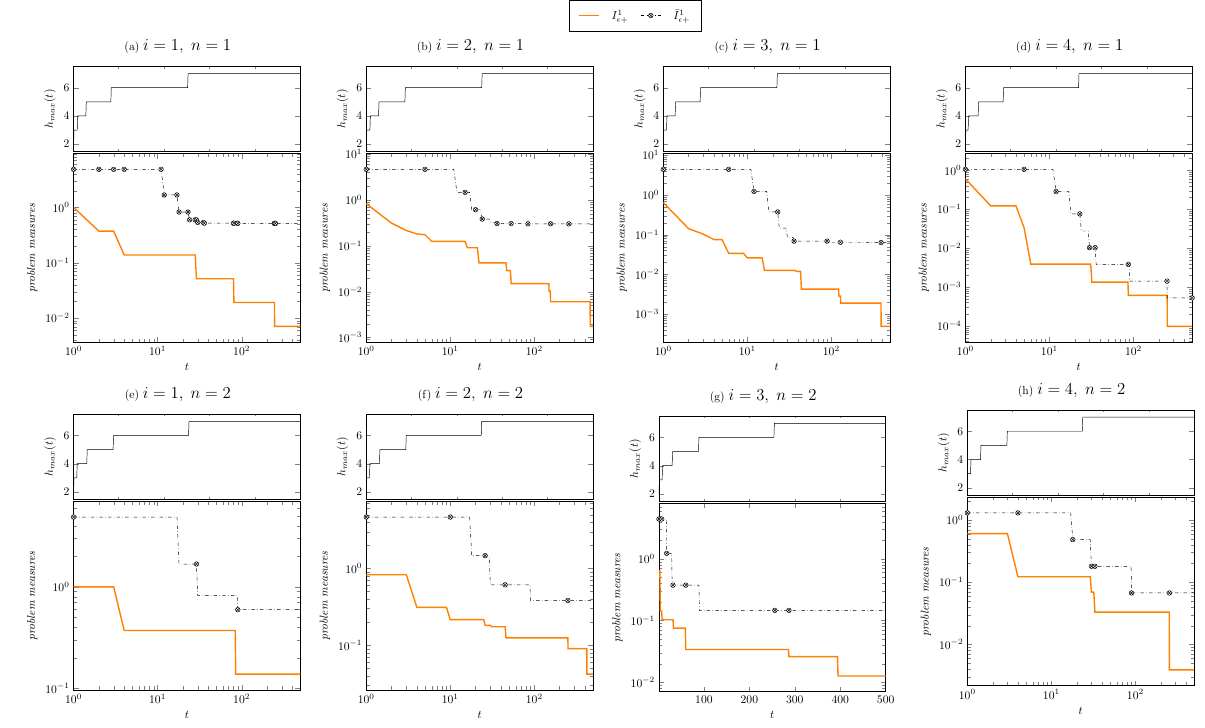}
		}
		\caption{Empirical Validation of Theorem~\ref{thm:finite-time-loss} using eight instances of problem~\eqref{eq:synthetic-problem} having the following configurations: (a) 
			$\vx^*_1=\mathbf{0},\;\vx^*_2=\mathbf{1}$, and $\mathcal{X}=[-1,1]^1$; 
			(b) $\vx^*_1=\mathbf{0.21}, \vx^*_2=\mathbf{0.81}$, and $\mathcal{X}=[-1,1]^1$; 
			(c) $\vx^*_1=\mathbf{0.47}, \vx^*_2=\mathbf{0.61}$, and $\mathcal{X}=[-1,1]^1$;
			(d) $\vx^*_1=\mathbf{0.57}, \vx^*_2=\mathbf{0.57}$, and $\mathcal{X}=[-1,1]^1$;
			(e) $\vx^*_1=\mathbf{0},\;\vx^*_2=\mathbf{1}$, and $\mathcal{X}=[-1,1]^2$; 
			(f) $\vx^*_1=\mathbf{0.21}, \vx^*_2=\mathbf{0.81}$, and $\mathcal{X}=[-1,1]^2$; 
			(g) $\vx^*_1=\mathbf{0.47}, \vx^*_2=\mathbf{0.61}$, and $\mathcal{X}=[-1,1]^2$;
			(h) $\vx^*_1=\mathbf{0.57}, \vx^*_2=\mathbf{0.57}$, and $\mathcal{X}=[-1,1]^2$. As a function of \ouralg's iterations, the figure shows three curves: the orange, solid curve denotes the additive $\epsilon$-indicator, while the dashed, dotted curve represents the quantity $\bar{I}^1_{\epsilon+}=\max_j \frac{1}{w_j}\delta(h(t), h_{max}(t) + 1)+ I^1_{\epsilon+}(\{\mathbf{f}(\arg\min_{\vx\in \mathcal{X}} \gfct(\vx))   \})$. In addition, the top region of each subplot shows the maximum tree depth reached $h_{\max}(t)$ at iteration $t$. The dashed curve upper bounds the orange, solid curve substantiating Theorem~\ref{thm:ind-bound}'s result: $ I^1_{\epsilon+}(\mathcal{Y}^t_*) - I^1_{\epsilon+}(\{\mathbf{f}(\arg\min_{\vx\in \mathcal{X}} \gfct(\vx))   \}) \leq \max_j \frac{1}{w_j}\delta(h(t), h_{max}(t) + 1)$.}
		\label{fig:emprcl_plts}
	\end{figure}
\end{landscape}

Figure~\ref{fig:emprcl_plts} presents the theoretical and numerical results obtained on the eight instances of problem~\eqref{eq:synthetic-problem}. It can be noted that the theoretical measures bound the numerical indicator throughout the algorithm iterations. Furthermore, the bound gets tighter (closer) on problems of less conflicting objectives ($i=3$ and $i=4$). That is, the optimal solutions of the objectives are close from each other in the decision space ---$(\mathbf{0.47}, \mathbf{0.61}), (\mathbf{0.57}, \mathbf{0.57})$, respectively. This is in line with the observation amount of exploration (the near-optimality dimension's constant ) grows linearly with the number of optimal solutions~\cite{Munos2014}. The code and data of the numerical validation will be made available at the project website.

\section{Conclusion}

This paper has established the Lipschitz continuity of the weighted \cheb function of Lipschitz-continuous multi-objective problems. The derived, yet unknown smoothness motivated formulating the weighted \cheb problem in a multi-armed bandits setting following the optimism in the face of uncertainty. As a result, we presented and described the Weighted Optimistic Optimization algorithm (\ouralg). \ouralg~looks for the optimal
solution by building deterministic hierarchical bandits (in the form of 
a space-partitioning tree) over the decision space $\mathcal{X}$. 

The presented finite-time analysis has established an upper bound on the Pareto-compliant additive $\epsilon$-indicator value of the approximation set obtained from \ouralg's sampled function values as a function of the number of iterations. Numerical experiments on a set of synthetic problems of varying difficulty and dimensionality confirmed the theoretical bounds. To the best of our knowledge, this is the first time, a Pareto-compliant quality indicator is investigated with respect to the weighted \cheb problem in finite time. Potential future research directions include studying the effect of randomizing the weighting vector $\mathbf{w}$ in a stochastic multi-armed bandits setting.

\bibliographystyle{spmpsci}
\bibliography{bib}

\end{document}